\documentclass[12pt]{article}
\usepackage{amsmath,amsfonts,amssymb,amsthm}
\begin{document}

\newcommand{\1}{{{\bf 1}}}
\newcommand{\id}{{\rm id}}
\newcommand{\Hom}{{\rm Hom}\,}
\newcommand{\End}{{\rm End}\,}
\newcommand{\Res}{{\rm Res}\,}
\newcommand{\Image}{{\rm Im}\,}
\newcommand{\Ind}{{\rm Ind}\,}
\newcommand{\Aut}{{\rm Aut}\,}
\newcommand{\Ker}{{\rm Ker}\,}
\newcommand{\gr}{{\rm gr}}
\newcommand{\Der}{{\rm Der}\,}
\newcommand{\Z}{\mathbb{Z}}
\newcommand{\Q}{\mathbb{Q}}
\newcommand{\C}{\mathbb{C}}
\newcommand{\N}{\mathbb{N}}
\newcommand{\g}{\mathfrak{g}}
\newcommand{\gl}{\mathfrak{gl}}
\newcommand{\h}{\mathfrak{h}}
\newcommand{\wt}{{\rm wt}\,}
\newcommand{\A}{\mathcal{A}}
\newcommand{\D}{\mathcal{D}}
\newcommand{\Lie}{\mathcal{L}}
\newcommand{\E}{\mathcal{E}}

\def \b{\beta}
\def \<{\langle}
\def \>{\rangle}
\def \be{\begin{equation}\label}
\def \ee{\end{equation}}
\def \bex{\begin{exa}\label}
\def \eex{\end{exa}}
\def \bl{\begin{lem}\label}
\def \el{\end{lem}}
\def \bt{\begin{thm}\label}
\def \et{\end{thm}}
\def \bp{\begin{prop}\label}
\def \ep{\end{prop}}
\def \br{\begin{rem}\label}
\def \er{\end{rem}}
\def \bc{\begin{coro}\label}
\def \ec{\end{coro}}
\def \bd{\begin{de}\label}
\def \ed{\end{de}}

\newtheorem{thm}{Theorem}[section]
\newtheorem{prop}[thm]{Proposition}
\newtheorem{coro}[thm]{Corollary}
\newtheorem{conj}[thm]{Conjecture}
\newtheorem{exa}[thm]{Example}
\newtheorem{lem}[thm]{Lemma}
\newtheorem{rem}[thm]{Remark}
\newtheorem{de}[thm]{Definition}
\newtheorem{hy}[thm]{Hypothesis}
\makeatletter \@addtoreset{equation}{section}
\def\theequation{\thesection.\arabic{equation}}
\makeatother \makeatletter

\begin{Large}
\begin{center}
\textbf{$C_2$-cofiniteness of the vertex algebra $V_L^+$ when $L$ is a non-degenerate even lattice}
\end{center}
\end{Large}
\begin{center}{
Phichet Jitjankarn\\ Department of Mathematics, Chulalongkorn University, Bangkok, Thailand,\\
and\\
Gaywalee Yamskulna\footnote{Address for correspondence: Department of Mathematical Sciences, Illinois State University, Normal, IL 61790, USA.\\
E-mail: gyamsku@ilstu.edu}\\
Department of Mathematical Sciences, Illinois State University, Normal, IL 61790\\
and\\
Institute of Science, Walailak University, Nakon Si Thammarat,
Thailand}
\end{center}
\begin{abstract}

In \cite{abd}, Abe, Buhl, and Dong showed that when $L$ is a positive definite even lattice, the vertex algebra $V_L^+$ and its irreducible weak modules satisfy the $C_2$-cofiniteness condition. In this paper, we extend their results by showing that the vertex algebra $V_L^+$ and its irreducible weak modules are $C_2$-cofinite when $L$ is a negative definite even lattice and when $L$ is a non-degenerate even lattice that is neither negative definite nor positive definite.\end{abstract}
\section{Introduction}

A vertex algebra $V$ is said to be $C_2$ cofinite if the subspace $C_2(V)=Span_{\C}\{~u_{-2}v~|~u,v\in V~\}$ has a finite co-dimension in $V$. This is often called the $C_2$ condition, and it was first appeared in \cite{z} when Zhu used this condition, as well as other assumptions, to show the modular invariance of certain trace functions. Since its introduction, the $C_2$-condition has proven to be a power tool in the study of theory of vertex algebras. In particular, it has played an important role 
in the study of structure of modules of vertex algebras which satisfy it (cf. \cite{bu, dlm3, gn, kl, m, m2}).

The vertex algebras $V_L^+$ are one of the most  important classes of vertex algebras along with those vertex algebras  associated with lattices, affine Lie algebras and Virasoro algebras. They were originally introduced in the Frenkel-Lepowsky-Meurman construction of the moonshine module vertex algebra (cf. \cite{flm}). When $L$ is a positive definite even lattice, the representation theory of the vertex algebra $V_L^+$ is well understood. In fact, for such case, the classification of all irreducible weak $V_L^+$-modules, and the study of the complete reducibility property of weak $V_L^+$-modules was done by Abe, Dong, Jiang, and Nagatomo (see \cite{a1, a, ad, dj, dn1}). Furthermore, the $C_2$-cofiniteness property of such vertex algebra $V_L^+$ and its irreducible weak modules had been shown by Abe, Buhl, Dong and the second author (cf. \cite{abd,y}). It turned out that this $C_2$-condition was one of the key ingredients for showing the complete reducibility of weak $V_L^+$-modules.

When $L$ is a non-degenerate even lattice that is not positive definite, the classification of weak $V_L^+$-modules are not completed yet. For the case when $L$ is a rank one negative definite even lattice, the classification of irreducible admissible $V_L^+$-modules was done by Jordan in \cite{j}. Later, in \cite{y2, y3}, the second author classified all irreducible admissible $V_L^+$-modules and showed the complete reducibility of admissible $V_L^+$-modules for the case when $L$ is a negative definite even lattice of arbitrary rank, and when $L$ is a non-degenerate even lattice that is neither positive definite nor negative definite. Note that an admissible module is a weak module that satisfies certain assumptions.

In this paper, we take a further step in understanding the representation theory of the vertex algebra $V_L^+$ when $L$ is a non-degenerate even lattice that is not positive definite by showing that such vertex algebra and its irreducible weak modules satisfy the $C_2$ condition. We hope to use these results in the study of the complete reducibility of weak $V_L^+$-modules when $L$ is a non-degenerate even lattice that is not positive definite in the future.

This paper is organized as follows. In Section 2, we review definitions of the vertex algebra, weak modules, and the cofiniteness $C_2$-condition. Also, we study vertex algebras that satisfy a special condition and the roles of this condition on weak modules of such vertex algebras. These results will play important roles in Section 4. In Section 3, we review the construction of the vertex algebras $V_L^+$ and their irreducible weak modules. Finally, in Section 4, we show that the vertex algebra $V_L^+$ and its irreducible weak modules satisfy the $C_2$-condition when $L$ is a negative definite even lattice and when $L$ is a non-degenerate even lattice that is neither positive definite nor negative definite. 

This work was done while the first author was visiting Illinois State University (ISU) from April 2008 to March 2009. P. J. thanks ISU Mathematics Department for their hospitality.

Throughout this paper, $\Z_{>0}$ is the set of positive integers.
\section{Vertex algebras and the $C_2$-condition}

First we review definitions of vertex algebras, weak modules and the $C_2$-condition. Next,  we show that when $V=C_2(V)$, many types of weak $V$-modules and the tensor-product-vertex algebra that contains $V$ as a vertex sub-algebra satisfy the $C_2$-condition. These will play important roles in Section 4.

\begin{de}\cite{lli} A {\em vertex algebra} $V$ is a vector space equipped with a linear map
$Y(\cdot, z):V\rightarrow(\End V)[[z,z^{-1}]], v\mapsto Y(v,z)=\sum_{n\in\Z}v_nz^{-n-1}$ and a distinguished vector ${\bf 1}\in V$ which satisfies the following properties: for $u,v\in V$
\begin{enumerate}
\item $u_nv=0$ for $n>>0$.
\item $Y({\bf 1},z)=id_V$.
\item $Y(v,z){\bf 1}\in V[[z]]$ and $\lim_{z\rightarrow 0}Y(v,z){\bf 1}=v$.
\item (the Jacobi identity) \begin{eqnarray*}
& &z_0^{-1}\delta\left(\frac{z_1-z_2}{z_0}\right)Y(u,z_1)Y(v,z_2)-z_0^{-1}\delta\left(\frac{z_2-z_1}{-z_0}\right)Y(v,z_2)Y(u,z_1)\\
&&=z_2^{-1}\delta\left(\frac{z_1-z_0}{z_2}\right)Y(Y(u,z_0)v,z_2).
\end{eqnarray*}
\end{enumerate}
We denote the vertex algebra just defined by $(V, Y, {\bf 1})$ or, briefly, by $V$.
\end{de}
\begin{de} A $\Z$-graded vertex algebra is a vertex algebra $$V=\oplus_{n\in\Z}V_n; \text{ for }v\in V_n,\ \  n=\wt v,$$ equipped with a conformal vector $\omega\in V_2$ which satisfies the following relations:
\begin{itemize}
\item $[L(m),L(n)]=(m-n)L(m+n)+\frac{1}{12}(m^3-m)\delta_{m+n,0}c_V$ for $m,n\in \Z$, where $c_V\in \C$ (the central charge) and $$Y(\omega,z)=\sum_{n\in\Z}L(n)z^{-n-2}\left(=\sum_{m\in\Z}\omega_mz^{-m-1}\right);$$
\item $L(0)v=nv=(\wt v) v$ for $n\in\Z$, and $v\in V_n$;
\item $Y(L(-1)v,z)=\frac{d}{dz}Y(v,z)$.
\end{itemize}
\end{de}
\begin{de}\cite{lli} A {\em vertex sub-algebra} of a vertex algebra $V$ is a vector sub-space $U$ of $V$  such that ${\bf 1}\in U$ and  such that $U$ is itself a vertex algebra.
\end{de}
\begin{de}\cite{lli} An {\em automorphism} of a vertex algebra $V$ is a linear isomorphism of $V$ such that $g({\bf 1})={\bf 1}$, and $gY(v,z)u=Y(g(v),z)g(u)$ for $u,v\in V$.
\end{de}
\begin{de}\cite{dlm2} A weak $V$-module $M$ is a vector space equipped with a linear map $Y_M(\cdot,z):V\rightarrow (\textnormal{End}M)[[z,z^{-1}]]$, $v\mapsto Y_M(v,z)=\sum_{n\in\Z}v_nz^{-n-1}$ which satisfies the following properties$:$ for $v,u\in V$, and $w\in M$
\begin{enumerate}
\item $u_nw=0$ for $n>>0$.
\item $Y({\bf 1},z)=id_M$.
\item (the Jacobi identity) \begin{eqnarray*}
& &z_0^{-1}\delta\left(\frac{z_1-z_2}{z_0}\right)Y_M(u,z_1)Y_M(v,z_2)-z_0^{-1}\delta\left(\frac{z_2-z_1}{-z_0}\right)Y_M(v,z_2)Y_M(u,z_1)\\
&&=z_2^{-1}\delta\left(\frac{z_1-z_0}{z_2}\right)Y_M(Y(u,z_0)v,z_2).
\end{eqnarray*}
\end{enumerate}
\end{de}

\begin{de}An {\em irreducible} weak $V$-module is a weak $V$-module that has no weak $V$-submodule except $\{0\}$ and itself. Here a weak submodule is defined in the obvious way.
\end{de}

\begin{de} A weak $V$-module is {\em completely reducible} if it can be rewritten as a direct sum of finitely many irreducible weak $V$-modules.
\end{de}

\begin{prop}\cite{lli} Let $W$ be a weak $V$-module, and let $\<T\>$ be a weak $V$-submodule of $W$ generated by a subset $T$ of $W$. Then $$\<T\>=Span_{\C}\{~v_nw~|~v\in V, n\in\Z, w\in T~\}.$$
\end{prop}
\begin{coro}\label{irrsp} If $W$ is an irreducible weak $V$-module then $$W=Span_{\C}\{~v_nw~|~v\in V, n\in\Z~\}.$$ Here, $w$ is a non-zero element in $W$.
\end{coro}





For the rest of this section, {\bf we assume that $V$ is a $\Z$-graded vertex algebra}.

\begin{de}\cite{z} We define $$C_2(V)=Span_{\C}\{~u_{-2}v~|~u,v\in V~\}.$$ 
$V$ is said to satisfy {\em the cofiniteness $C_2$-condition} if $V/C_2(V)$ has a finite dimension. \end{de}

\begin{prop}\label{prop2.7}\cite{z} 
\begin{enumerate}
\item $L(-1)u\in C_2(V)$ for $u\in V$.
\item $v_{-n}u\in C_2(V)$ for $u,v\in V$ and $n\geq 2$.
\item $C_2(V)$ is an ideal of $V$ with respect to $(-1)^{st}$-product.
\item $V/C_2(V)$ is a commutative associative algebra under $(-1)^{st}$-product.
\end{enumerate}
\end{prop}
\begin{de} Let $M$ be a weak $V$-module. We define
$$\widetilde{C}_2(M)=Span_{\C}\{~u_{-2}w~|~u\in V,~w\in M~\}.$$
The space $M$ is called \em{$C_2$-cofinite} if $M/\widetilde{C}_2(M)$ has a finite dimension. 
\end{de}
The following is the key proposition.
\begin{prop}\label{prop4.3} Let $W$ be an irreducible weak $V$-module. If $V=C_2(V)$, then $W=\widetilde{C}_2(W)$.
\end{prop}
\begin{proof} We first show that $a_nw\in \widetilde{C}_2(W)$ for all $a\in V, w\in W$ and $n\in\Z$. Clearly, this statement is true when $n\leq -2$. Since $V=C_2(V)$, one may write $a_nw$ as
$$a_{n}w=\sum_{j=1}^l(u^j_{-2}v^j)_nw=\sum_{j=1}^l\sum_{i\geq0}(-1)^i \binom{-2}{i} \left(u^j_{-2-i}v^j_{n+i}-v^j_{-2+n-i}u^j_i\right)w.$$ Here $u^j,v^j\in V$. By using an induction on $n$, we can show that $a_nw\in \widetilde{C}_2(W)$ for all $a\in V, w\in W$ and $n\geq -1$. 

Next, by using the facts that $W$ is irreducible and $a_nw\in \widetilde{C}_2(W)$ for all $a\in V, w\in W$, and $n\in\Z$, we can conclude immediately that $W=\widetilde{C}_2(W)$.
\end{proof}
\begin{coro}\label{coro2.13} Assume that $V=C_2(V)$. If a weak $V$-module $M$ is completely reducible, then $M=\widetilde{C}_2(M)$.
\end{coro}
\begin{prop}\label{prop2.14} Let $V'$ be a $\Z$-graded vertex algebra and let $V$ be a $\Z$-graded vertex sub-algebra of $V'$ such that $V=C_2(V)$. If $V'$ is completely reducible as a weak $V$-module, then $V'$ is $C_2$-cofinite. In particular, $V'=C_2(V')$.
\end{prop}
\begin{proof} By Corollary \ref{coro2.13}, we can conclude that $V'=\widetilde{C}_2(V')$. Since $\widetilde{C}_2(V')$ is a subset of $C_2(V')$, we then have that $V'=C_2(V')$.
\end{proof}
\begin{prop} Let $V^1,...,V^n$ be $\Z$-graded vertex algebras. Assume that $V^j=C_2(V^j)$ for some $j\in \{1,\cdots, n\}$. Then $V^1\otimes...\otimes V^n$ satisfies the $C_2$-condition. In particular, $V^1\otimes...\otimes V^n=C_2(V^1\otimes...\otimes V^n)$.
\end{prop}
\begin{proof} The statements follow from the fact that 
\begin{eqnarray*}
V^1\otimes...\otimes V^n&=&V^1\otimes .\cdots\otimes C_2(V^j)\otimes\cdots\otimes V^n\\&\subset& C_2(V^1\otimes...\otimes V^n).
\end{eqnarray*}
\end{proof}

\section{The vertex algebras $V_L^+$ and their irreducible weak modules}

In this section, we review the construction of the vertex algebra $V_L^+$ and its irreducible weak $V_L^+$-modules.

Let $L$ be a rank $d$ even lattice equipped with a non-degenerate symmetric $\Z$-bilinear form $\<\cdot,\cdot \>$. Set $\mathfrak{h}=\C\otimes_{\Z}L$ and extend $\<\cdot,\cdot\>$ to a $\C$-bilinear form on $\h$. Let $\mathfrak{\widehat{h}}=\mathfrak{h}\otimes \C [t,t^{-1}] \oplus \C c$ be the corresponding affine Lie algebra 
with the following commutator relations:\[\left[\beta\otimes t^m,\gamma\otimes t^n\right]=\left\langle \beta,\gamma\right\rangle m \delta_{m+n,0}c~~\text{and}~~[c,\mathfrak{\widehat{h}}]=0\]
where $\beta,\gamma\in\mathfrak{h}$ and $m,n\in\Z$. We consider the induced $\hat{\h}$-module
\[M(1) = U(\mathfrak{\widehat{h}})\otimes_{U(\mathfrak{h}\otimes\C[t]\oplus\C)} \C. \]
Here, $\mathfrak{h}\otimes \C[t]$ acts trivially on $\C$, and  $c$ acts as a multiplication by 1. Then there exists a linear map $Y:M(1)\rightarrow\End(M(1))[[z,z^{-1}]]$ such that $M(1)$ becomes a simple $\Z$-graded vertex algebra with the vacuum vector $1$ and the Virasoro element $\omega=\frac{1}{2}\sum^d_{i=1}\beta_i(-1)^2\otimes1$ (see \cite{flm}). Here, $\left\{\beta_1,\ldots,\beta_d\right\}$ is an orthonormal basis of $\h$.


Next, we let $\hat{L}$ be a canonical central extension of $L$ by the cyclic group of order 2:$$1\rightarrow\<\kappa\>\rightarrow \hat{L}\rightarrow L\rightarrow 0$$ with the commutator map $c(\alpha,\beta)=\kappa^{\<\alpha,\beta\>}$ for $\alpha,\beta\in L$. Let $e:L\rightarrow\hat{L}$ be a section such that $e_0=1$ and let $\epsilon:L\times L\rightarrow \<\kappa\>$ be the corresponding bimultiplicative 2-cocycle. Then $\epsilon(\alpha,\beta)\epsilon(\beta,\alpha)=\kappa^{\<\alpha,\beta\>}$, $\epsilon(\alpha,\beta)\epsilon(\alpha+\beta,\gamma)=\epsilon(\beta,\gamma)\epsilon(\alpha,\beta+\gamma)$ and $e_{\alpha}e_{\beta}=\epsilon(\alpha,\beta)e_{\alpha+\beta}$ for $\alpha,\beta,\gamma\in L$.

Let $L^{\circ}=\{~\lambda\in\h~|~\<\alpha,\lambda\>\in \Z\text{ for all }\alpha\in L~\}$ be a dual lattice of $L$. There is a $\hat{L}$-module structure on $\C[L^{\circ}]=\oplus_{\lambda\in L^{\circ}}\C e^{\lambda}$ such that $\kappa$ acts as $-1$. We set $L^{\circ}=\cup_{i\in L^{\circ}/L}(L+\lambda_i)$ be the coset decomposition such that $\lambda_0=0$, and we define
$\C[L+\lambda_i]=\oplus_{\alpha\in L}\C e^{\alpha+\lambda_i}$. Clearly, $\C[L+\lambda_i]$ is a ${L^{\circ}}$-module. Moreover, $\C[L+\lambda_i]$ is a $\hat{L}$-module under the following action: $e_{\alpha}e^{\beta+\lambda_i}=\epsilon(\alpha,\beta)e^{\alpha+\beta+\lambda_i}$ for $\alpha,\beta\in L$.

For simplicity, we will identify $e^{\alpha}$ with $e_{\alpha}$ for $\alpha\in L$. Also, for a subset $M$ of $L^{\circ}$, we set $\C[M]=\oplus_{\lambda\in M}\C e^{\lambda}$. Let $z$ be a formal variable and $h\in\h$. For $\lambda\in L^{\circ}$, we define actions of $h$ and $z^h$ on $\C[L+\lambda]$ in the following ways: $h\cdot e^{\beta}=\<h,\beta\>e^{\beta};$ and $z^h\cdot e^{\beta}=z^{\<h,\beta\>}e^{\beta}$. Now, we set $$V_{L+\lambda}=M(1)\otimes \C[L+\lambda].$$ Clearly $L$, $\hat{\h}$, $z^h$ ($h\in\h$) will act naturally on $V_{L+\lambda}$ by acting either on $M(1)$ or $\C[L+\lambda]$.  
For $h\in \h$, $\beta\in L$, the vertex operators $Y(h(-1){\bf 1},z)$ and $Y(e^{\beta},z)$ associated to $h(-1){\bf 1}$ and $e^{\beta}$, respectively, are defined as follow:
\begin{eqnarray*}
& &Y(h(-1){\bf 1},z)=h(z)=\sum_{n\in\Z}h(n)z^{-n-1}\\
& &Y(e^{\beta},z)=\exp\left(\sum_{n=1}^{\infty}\frac{\beta(-n)}{n}z^n\right)\exp\left(-\sum_{n=1}^{\infty}\frac{\beta(n)}{n}z^{-n}\right)e^{\beta}z^{\beta}.
\end{eqnarray*}
The vertex operator associated with a vector $v=\gamma_1(-n_1)...\gamma_r(-n_r)e^{\beta}$ for $\gamma_i\in \h$, $n_i\geq 1$ and $\beta\in L$ is defined by
$$Y(v,z)=:\partial^{(n_1-1)}\gamma_1(z)...\partial^{(n_r-1)}\gamma_r(z)Y(e^{\beta},z):,$$
where $\partial^{(n)}=\frac{1}{n!}\left(\frac{d}{dz}\right)^n$ and the normal ordering $:\cdot :$ is an operation which reorders the operators so that $\gamma(n)$ ($\gamma\in\h,n<0$) and $e^{\beta}$ to be placed to the left of $\gamma(n)$ ($\gamma\in\h, n\geq 0$) and $z^{\beta}$.

\begin{thm}\cite{b,d1,flm} 

\begin{enumerate}
\item The space $V_L$ is a simple $\Z$-graded vertex algebra with a Virasoro element $\omega=\frac{1}{2}\sum^d_{i=1}h_i(-1)^2\otimes e^0$. Here $\left\{h_i~|~1\leq i\leq d\right\}$ is an orthonormal basis of $\h$. Moreover, $M(1)$ is a $\Z$-graded vertex sub-algebra of $V_L$. 
\item $\{~V_{L+\lambda_i}~|~i\in L^{\circ}/L~\}$ is the set of all inequivalent irreducible $V_L$-modules.
\end{enumerate}
\end{thm}

Next, we define a linear isomorphism $\theta:V_{L+\lambda_i}\rightarrow V_{L-\lambda_i}$ by
$$\theta(\b_1(-n_1)\b_2(-n_2)....\b_k(-n_k)e^{\beta+\lambda_i})=(-1)^k\b_1(-n_1)....\b_k(-n_k)e^{-\beta-\lambda_i}$$ if $2\lambda_i\not\in L$ and 
\begin{eqnarray*}
& &\theta(\b_1(-n_1)\b_2(-n_2)....\b_k(-n_k)e^{\beta+\lambda_i})\\
& &=(-1)^kc_{2\lambda_i}\epsilon(\beta,2\lambda_i)\b_1(-n_1)....\b_k(-n_k)e^{-\beta-\lambda_i}
\end{eqnarray*} if $2\lambda_i\in L$. Here $c_{2\lambda_i}$ is a square root of $\epsilon(2\lambda_i,2\lambda_i)$. 
Consequently, $\theta$ is a linear automorphism of $V_{L^{\circ}}$. Moreover, $\theta$ is an automorphism of $V_L$ and $M(1)$.

For any stable $\theta$-subspace $U$ of $V_{L^{\circ}}$, we denote a $\pm 1$ eigenspace of $U$ for $\theta$ by $U^{\pm}$.
\begin{prop}\cite{dm, dlm1} 

\begin{enumerate}
\item $M(1)^+$ and $V_L^+$ are simple $\Z$-graded vertex algebras.
\item $(V_{L+\lambda_i}+V_{L-\lambda_i})^{\pm}$ for $i\in L^{\circ}/L$ are irreducible weak $V_L^+$-modules. Moreover, if $2\lambda_i\not\in L$ then $(V_{L+\lambda_i}+V_{L-\lambda_i})^{\pm}$, $V_{L+\lambda_i}$ and $V_{L-\lambda_i}$ are isomorphic.
\end{enumerate}
\end{prop}

Let $\h[-1]=\h\otimes t^{\frac{1}{2}}\C[t,t^{-1}]\oplus \C c$ be the twisted affine Lie algebra with the commutator relation
\[\left[\beta\otimes t^m,\gamma\otimes t^n\right]=\left\langle \beta,\gamma\right\rangle m \delta_{m+n,0}c~~\text{and}~~[c,\h[-1]]=0\]
where $\beta,\gamma\in\mathfrak{h}$ and $m,n\in\Z+\frac{1}{2}$. Next we let 
$M(1)(\theta)=S(\h\otimes t^{-\frac{1}{2}}\C[t^{-1}])$ be the unique irreducible $\h[-1]$-module such that $c$ acts as 1, and when $n>0$, $\beta\otimes t^n$ acts on 1 as zero.

By abusing the notation, we use $\theta$ to denote an automorphism of $\hat{L}$ defined by $\theta(e_\alpha)=e_{-\alpha}$ and $\theta(\kappa)=\kappa$. We set $K=\left\{a^{-1}\theta(a)~|~a\in\hat{L}\right\}$. Also, we let $\chi$ be a central character of $\hat{L}/K$ such that $\chi(\kappa)=-1$ and we let $T_\chi$ be the irreducible $\hat{L}/K$-module with central character $\chi$. We define $$V^{T_\chi}_L=M(1)(\theta)\otimes T_\chi.$$ 
We define an action of $\theta$ on $V^{T_\chi}_L$ in the following way:
\begin{center}
$\begin{array}{lcl}
\theta(\b_1(-n_1)\b_2(-n_2)....\b_k(-n_k)t) &=& (-1)^k\b_1(-n_1)....\b_k(-n_k)t.    
\end{array}$
\end{center}
Here, $\beta_i\in\h$, $n_i\in \frac{1}{2}+\Z_{\geq0}$ and $t\in T_\chi$. We denote by $V^{T_\chi,\pm}_L$ the $\pm1$-eigenspace of $V^{T_\chi}_L$ for $\theta$.

\begin{prop}\cite{dli} Let $\chi$ be a central character of $\hat{L}/K$ such that $\chi(\iota(\kappa))=-1$ and let $T_{\chi}$ be the irreducible $\hat{L}/K$-module with central character $\chi$. Then $V_L^{T_{\chi},\pm}$ are irreducible weak $V_L^+$-modules. \end{prop}


\section{The $C_2$-cofiniteness condition of the vertex algebra $V_L^+$ and its weak modules}

When $L$ is a positive definite even lattice, it was shown by Abe, Buhl and Dong, and the second author that the vertex algebras $V_L^+$ and their irreducible weak modules satisfy the $C_2$ condition (cf. \cite{abd,y}). In this section, we will extend their results to the case when $L$ is a negative definite even lattice and when $L$ is a non-degenerate even lattice that is neither positive definite nor negative definite.

\subsection{Case I: $L$ is a rank one negative definite even lattice.}
For the rest of this subsection, we assume that $L=\Z\alpha$ is a rank one negative definite even lattice such that $\<\alpha,\alpha \>=-2k$. Here, $k$ is a positive integer.

Let $m\in\Z_{>0}$. For convenience, we set 
\begin{eqnarray*}
& & E^m=e^{m\alpha}+e^{-m\alpha}, \ \ F^m=e^{m\alpha}-e^{-m\alpha},\text{ and }\\
& &V_L^+(m)=M(1)^+\otimes E^m\oplus M(1)^-\otimes F^m.\end{eqnarray*}
Clearly, $V_{L}^+=M(1)^+\oplus \oplus_{m=1}^{\infty}V_L^+(m).$ 

\begin{prop}\label{prop2.11}\cite{dg, dn1, j} \ \ 
\begin{enumerate}
\item As a vertex algebra, $M(1)^+$ is generated by the Virasoro element $\omega$ and any singular vector of weight greater than zero. In particular, $M(1)^+$ is generated by $\omega$ and $J$ where $J=\frac{1}{4k^2}\alpha(-1)^4\mathbf{1}+\frac{1}{k}\alpha(-3)\alpha(-1)\mathbf{1}-\frac{3}{4k}\alpha(-2)^2\mathbf{1}$.

\item The vertex algebra $M(1)^+$ is spanned by $$L(-m_1)\ldots L(-m_s)J_{-n_1}\ldots J_{-n_t}\mathbf{1}$$ 
where $m_1\geq\cdots\geq m_s\geq2$, $n_1\geq\cdots\geq n_t\geq1$.

\item $V_L^+(m)$ is spanned by $L(-n_1)L(-n_2)....L(-n_r)\otimes E^{m}$ where $n_r\geq ....\geq n_2\geq n_1\geq 1$.
 
\end{enumerate}
\end{prop}

Next, for $m>0$, we set 
$$V_L^+(m)=\oplus_{n\in\Z} V_L^+(m,n)$$
where $V^+_L(m,n)$ is the weight $n$ subspace of $V^+_L(m)$. 
Clearly, $V_L^+(m,-km^2+6)$ has the following basis elements:

\underline{Basis of $V^+_L(m,-km^2+6)$ }

\[\begin{array}{lll}
g_1=\alpha(-6)F^m,&g_2=\alpha(-5)\alpha(-1)E^m,&g_3=\alpha(-4)\alpha(-2)E^m,\\              g_4=\alpha(-4)\alpha(-1)^2F^m,&g_5=\alpha(-3)^2E^m,&g_6=\alpha(-3)\alpha(-2)\alpha(-1)F^m,\\
g_7=\alpha(-3)\alpha(-1)^3E^m,&g_8=\alpha(-2)^3F^m,&g_9=\alpha(-2)^2\alpha(-1)^2E^m,\\
g_{10}=\alpha(-2)\alpha(-1)^4F^m,&g_{11}=\alpha(-1)^6E^m.
\end{array}\]
Furthermore, the following elements form bases of $V^+_L(m,-km^2+5)$ and $V^+_L(m,-km^2+3)$, respectively:

\ \ 

\underline{Basis of $V^+_L(m,-km^2+5)$ }
\ \ 
\begin{center}
$\begin{array}{lll}
f_1=\alpha(-5)F^m,&f_2=\alpha(-4)\alpha(-1)E^m,&f_3=\alpha(-3)\alpha(-2)E^m,\\
f_4=\alpha(-3)\alpha(-1)^2F^m,&f_5=\alpha(-2)^2\alpha(-1)F^m,&f_6=\alpha(-2)\alpha(-1)^3E^m,\\
f_{7}=\alpha(-1)^5F^m
\end{array}$
\end{center}

\underline{Basis of $V^+_L(m,-km^2+3)$}
\begin{center}
$\begin{array}{lll}
h_1=\alpha(-3)F^m,&h_2=\alpha(-2)\alpha(-1)E^m,&h_3=\alpha(-1)^3F^m.
\end{array}$
\end{center}

\begin{lem}\label{l3.1} The set $\{~L(-1)f_i,~L(-3)h_j, ~(\alpha(-1)^4\mathbf{1})_{-3}E^m~|~ 1\leq i\leq 7, ~ 1\leq j\leq 3\}$ is a basis of $V^+_L(m,-km^2+6).$
\end{lem}
\begin{proof} The table 1 describes expressions of $L(-1)f_i$~~$(i=1,\ldots,7)$, $L(-3)h_j$~~$(j=1,2,3)$, and $(\alpha(-1)^4\mathbf{1})_{-3}E^m$ in terms of $g_l$ $(l=1,\ldots,11)$. We will denote this table by a $11\times 11$ matrix $A$. Since $\det A=-6144m^3k^2\left(m^2k+1\right)$, and $m, k$ are positive integers, we can conclude that $A$ is invertible. Moreover, $L(-1)f_i$~~$(i=1,\ldots,7)$, $L(-3)h_j$~~$(j=1,2,3)$, and $(\alpha(-1)^4\mathbf{1})_{-3}E^m$ form a basis of $V^+_L(m,-km^2+6).$\end{proof}
  \begin{center}
  {\footnotesize
  \begin{tabular}{| c | @{}c@{} | @{}c@{} | @{}c@{} | @{}c@{} | @{}c@{} | @{}c@{} | @{}c@{} | @{}c@{} | @{}c@{} | @{}c@{} | @{}c@{} |}
    \hline
             & $g_1$ & $g_2$ & $g_3$ & $g_4$ & $g_5$ & $g_6$ & $g_7$ & $g_8$ & $g_9$ & $g_{10}$ & $g_{11}$ \\ \hline
    $L(-1)f_1$ & 5   & $m$   &  0  &  0  &  0  &  0  &  0  &  0  &  0  &  0     &  0      \\ 
    $L(-1)f_2$ & 0   & 4   &  1  &  $m$  &  0  &  0  &  0  &  0  &  0  &  0     &  0      \\ 
    $L(-1)f_3$ & 0   & 0   &  3  &  0  &  2  &  $m$  &  0  &  0  &  0  &  0     &  0      \\
    $L(-1)f_4$ & 0   & 0   &  0  &  3  &  0  &  2  &  $m$  &  0  &  0  &  0     &  0      \\
    $L(-1)f_5$ & 0   & 0   &  0  &  0  &  0  &  4  &  0  &  1  &  $m$  &  0     &  0      \\
    $L(-1)f_6$ & 0   & 0   &  0  &  0  &  0  &  0  &  2  &  0  &  3  &  $m$     &  0      \\
    $L(-1)f_7$ & 0   & 0   &  0  &  0  &  0  &  0  &  0  &  0  &  0  &  5     &  $m$      \\
    $2kL(-3)h_1$ & $6k$   & 0   &  0  &  0  &  $2km$  &  -1  &  0  &  0  &  0  &  0     &  0      \\
    $2kL(-3)h_2$ & 0   & $4k$   &  $2k$  &  0  &  0  &  $2km$  &  0  &  0  &  -1  &  0     &  0      \\
    $2kL(-3)h_3$ & 0   & 0   &  0  &  $6k$  &  0  &  0  &  $2km$  &  0  &  0  &  -1     &  0      \\
    $(\alpha(-1)^4\mathbf{1})_{-3}E^m$ & $-32k^3m^3$   & $48k^2m^2$   &  $48k^2m^2$  &  $-24km$  &  $24k^2m^2$  &  $-48km$  &  4  &  $-8km$  &  6  &  0     &  0      \\\hline
  \end{tabular}\smallskip\\Table 1.}
  \end{center}\smallskip

 \begin{prop}\label{prop3.2} For $m\in\Z_{>0}$, we have $V^+_L(m,-km^2+2n)$ is a subset of $C_2(V^+_L)$ when $n\geq 3$.
 \end{prop}
\begin{proof} By Proposition \ref{prop2.7} and Lemma \ref{l3.1}, we can conclude immediately that $V_L^+(m,-km^2+6)$ is contained in $C_2(V_L^+)$. 

Next, we will show that $V^+_L(m,-km^2+2n)$ is a subset of $C_2(V_L^+)$ for $n\geq 4$.
Since $(L(-2))^3E^m\in C_2(V_L^+)$, it implies that for $j\geq 4$, $(L(-2))^j E^m=(L(-2))^{j-3}(L(-2))^3E^m\in C_2(V_L^+)$. By Proposition \ref{prop2.11}, we can conclude further that $V^+_L(m,-km^2+2n)$ is a subset of $C_2(V^+_L)$ when $n\geq 4$.
\end{proof}
\begin{thm} For $m\in Z_{>0}$, $V_L^+(m)$ is a subset of $C_2(V_L^+)$.\end{thm}
\begin{proof} First, we set $\exp(\sum_{n=1}^{\infty}\frac{x_n}{n}z^{n}) = \sum_{j=0}^{\infty}p_j(x_1,x_2,...)z^{j}.$ Since  
\begin{eqnarray*} 
E_{-2kn}\alpha(-1)F^n&=&p_{4nk-1}(\alpha)\alpha(-1)\otimes e^{\alpha(n+1)}-p_{4nk-1}(-\alpha)\alpha(-1)\otimes e^{-\alpha(n+1)}\\
& & +2k\left(p_{4nk}(\alpha)\otimes e^{\alpha(n+1)}+p_{4nk}(-\alpha)\otimes e^{-\alpha(n+1)}\right)-2kE^{n-1}, 
\end{eqnarray*} 
and for $n\geq 2$, $V^+_L(n+1,-k(n+1)^2+4nk)$ is contained in $C_2(V^+_L)$ (cf. Proposition \ref{prop3.2}), it follows that $E^{n-1}\in C_2(V_L^+)$ for $n\geq 2$. By Proposition \ref{prop2.11} we can conclude further that $V_L^+(m)$ is a subset of $C_2(V_L^+)$ for all $m\in\Z_{>0}$.
\end{proof}

\begin{thm}\label{thm3.6} $V_L^+$ satisfies the $C_2$-cofiniteness condition. In particular, $$V_L^+/C_2(V_L^+)=\{0+C_2(V_L^+)\}.$$
\end{thm}
\begin{proof} We will show that $V_L^+=C_2(V_L^+)$. It is enough to show that $M(1)^+$ is contained in $C_2(V_L^+)$. Since $E_{-2k-1}E=p_{4k}(\alpha)\otimes e^{2\alpha}+p_{4k}(-\alpha)\otimes e^{-2\alpha}+2$ and $V_L^+(2)$ is a subset of $C_2(V_L^+)$, these imply that ${\bf 1}\in C_2(V_L^+)$. By Proposition \ref{prop2.11}, we can conclude that $M(1)^+$ is a subset of $C_2(V_L^+)$. Consequently, $V_L^+=C_2(V_L^+)$.\end{proof}
\begin{coro} Every irreducible weak $V_L^+$-module satisfies the $C_2$-condition. 
\end{coro}
\begin{proof} This follows immediately from Proposition \ref{prop4.3} and Theorem \ref{thm3.6}. \end{proof}

\subsection{Case II: $L$ is a negative definite even lattice}

For the rest of this subsection, we assume that $L$ is a rank $d$ negative definite even lattice.
\begin{thm}\label{th3.7} The vertex algebra $V_L^+$ is $C_2$-cofinite.\end{thm}
\begin{proof} We will show that $V_L^+=C_2(V_L^+)$. We will follow the proof in \cite{abd} very closely.
Let $K$ be a direct sum of $d$ orthogonal rank one negative definite even lattice $L_1,\ldots,L_d$. We write $L=\cup_{i\in L/K}(K+\lambda_i)$ as a direct sum of its coset decompositon with respect to $K$. Then $V_L=\oplus_{i\in L/K}V_{K+\lambda_i}$.

Since
\[V_K=V_{L_1}\otimes\cdots\otimes V_{L_d}=\displaystyle\sum_{\epsilon_i=\pm}V_{L_1}^{\epsilon_1}\otimes\cdots\otimes V_{L_d}^{\epsilon_d},\] 
it follows that
$V_K^+=\displaystyle\sum_{\epsilon_i=\pm,~\prod_i|\epsilon_i|=1}V_{L_1}^{\epsilon_1}\otimes\cdots\otimes V_{L_d}^{\epsilon_d}.$  Here $|\pm|=\pm1$. 
By using the fact that $V_{L_j}^+=C_2(V_{L_j}^+)$ for all $1\leq j\leq d$, we then have $$V_{L_1}^+\otimes\cdots\otimes V_{L_d}^+=C_2\left(V_{L_1}^+\otimes\cdots\otimes V_{L_d}^+\right).$$ Since $V_{L_1}^{\epsilon_1}\otimes\cdots\otimes V_{L_d}^{\epsilon_d}$ is an irreducible weak $V_{L_1}^+\otimes\cdots\otimes V_{L_d}^+$-module, it follows from Proposition \ref{prop2.14} that $V_K^+
={C}_2(V_K^+)$. Since $V_L^+$ is completely reducible as a weak $V_K^+$-module, we can conclude that $V_L^+=C_2(V_L^+)$ (cf. Proposition \ref{prop2.14}).
\end{proof}

\begin{coro} Every irreducible weak $V_L^+$-module satisfies the $C_2$-condition.\end{coro}

\subsection{Case III: $L$ is a rank $d$ non-degenerate even lattice that is neither positive definite nor negative definite}

\begin{thm} Assume that $L$ is a non-degenerate even lattice of a finite rank that is neither positive definite nor negative definite. Then the vertex algebra $V_L^+$ and its irreducible weak modules satisfy the $C_2$-condition. In particular, we have $V_L^+=C_2(V_L^+)$ and $M=\widetilde{C}_2(M)$ for all irreducible weak $V_L^+$-module.
\end{thm}
\begin{proof} The proof is very similar to the proof of Theorem \ref{th3.7}. In fact, we can obtain the above results by using the fact that $L$ has a rank one negative definite even sublattice, and follow the proof in Theorem \ref{th3.7} step by step.
\end{proof}


\begin{thebibliography}{BFFLS}
\bibitem[A1]{a1} T. Abe, The charge conjugation orbifold $V_{\Z\alpha}^+$ is rational when $\<\alpha, \alpha\>/2$ is prime, {\em Internat. Math. Res. Notices} {\bf 12} (2002), 647Ð665.
\bibitem[A2]{a} T. Abe, Rationality of the vertex operator algebra $V_L^+$ for a positive definite even lattice $L$, {\em Math. Z.} {\bf 249} (2005), no. 2, 455-484.
\bibitem[ABD]{abd}T. Abe, G. Buhl, and C. Dong, Rationality, regularity, and $C_2$-cofiniteness, {\em Trans. Amer. Math. Soc.}  {\bf 356} (2004) 3391-3402.
\bibitem[AD]{ad} T. Abe and C. Dong, Classification of irreducible modules for the vertex operator algebra $V_L^+$: General case, {\em J. Algebra} {\bf 273} (2004) no. 2, 657-685.
\bibitem[B]{b} R. Borcherds, Vertex algebras, Kac-Moody algebras and the Monster, {\em Proc. Natl. Acad. Sci. USA} {\bf 83} (1986), 3068-3071.
\bibitem[Bu]{bu} G. Buhl, A spanning set for VOA modules, {\em J. Algebra} {\bf 254} (2002) no. 1, 125-151.
\bibitem[D]{d1} C. Dong, Vertex algebras associated with even lattices, {\em J. Algebra} {\bf 160} (1993), 245-265.
\bibitem[DG]{dg} C. Dong and R. Griess, Rank one lattice type vertex operator algebras and their automorphism groups, {\em J. Algebra} {\bf 208} (1998), no. 1, 262-275.
\bibitem[DJ]{dj}
C. Dong and C. Jiang, Rationality of vertex operator algebras, math.QA/0607679.
\bibitem[DLM1]{dlm1} C. Dong, H.-S. Li and G. Mason, Compact automorphism groups of vertex operator algebras, {\em International Math. Research Notices} {\bf 18} (1996), 913-921.
\bibitem[DLM2]{dlm2} C. Dong, H.-S. Li and G. Mason, Regularity of rational vertex operator algebras, {\em Adv. Math.} {\bf 132} (1997), no. 1, 148--166.
\bibitem[DLM3]{dlm3} C. Dong, H.-S. Li and G. Mason, Modular invariance of trace functions in orbifold theory, {\em CMP} {\bf 214}, 1-56.
\bibitem[DLi]{dli} C. Dong and Z. Lin, Induced modules for vertex operator algebras, {\em Comm. Math. Phys.} {\bf 179} (1996), 157-184.
\bibitem[DM]{dm}
C. Dong and G. Mason, On quantum Galois theory, {\em Duke Math. J.} {\bf 86} (1997), 305-321.
\bibitem[DN]{dn1} 
C. Dong and K. Nagatomo, Representations of Vertex operator algebra $V_L^+$ for a rank one lattice $L$, {\em Comm. Math. Phys.} {\bf 202} (1999), 169-195.
\bibitem[FLM]{flm}
I.B. Frenkel, J. Lepowsky and A. Meurman, Vertex Operator Algebras and the Monster, Pure and Applied Math. Vol. {\bf 134}, Academic Press, Boston, 1988.
\bibitem[GN]{gn} M. Gaberdiel and A. Neitzke, Rationality; quasi rationality, and finite $W$-algebras, {\em Commun. Math. Phys.} {\bf B168}, 407-436.
\bibitem[J]{j}
L. Jordan, Classification of irreducible $V_L^+$-modules for a negative definite rank one even lattice $L$, Ph. D. Dissertation, University of California at Santa Cruz, 2006.
\bibitem[KL]{kl} M. Karel and H.-S. Li, Certain generating subspace for vertex operator algebras, {\em J. Algebra} {\bf 217}, 495-514.
\bibitem[LLi]{lli} J. Lepowsky and H.-S. Li, Introduction to Vertex Operator Algebras and Their Representations, Progress in Math. {\bf 227},
Birkh\"auser, Boston, 2003.
\bibitem[M1]{m} M. Miyamoto, Modular invariance of vertex operator algebras satisfying $C_2$-cofiniteness, {\em Duke Math. J.} {\bf 122} (2004), no. 1, 51-91.
\bibitem[M2]{m2} M. Miyamoto, A theory of tensor products for vertex operator algebra satisfying $C_2$-cofiniteness, math.QA/0309350.
\bibitem[Y1]{y} G. Yamskulna, $C_2$-cofiniteness of vertex operator algebra $V_L^+$ when $L$ is a rank one lattice. {\em Comm. Algebra.} {\bf 32} (2004), 927-954.
\bibitem[Y2]{y2} G. Yamskulna, Classification of irreducible modules of the vertex algebra $V_L^+$ when $L$ is a nondegenerate even lattice of an arbitrary rank, {\em J. Algebra} {\bf 320} (2008), 2455-2480.
\bibitem[Y3]{y3} G. Yamskulna, Rationality of the vertex algebra $V_L^+$ when $L$ is a non-degenerate even lattice of arbitrary rank, {\em J. Algebra} {\bf 321} (2009), 1005-1015.
\bibitem[Z]{z} Y. Zhu, Modular invariance of characters of vertexoperator algebras, {\em J. Amer. Math. Soc.} {\bf 9} (1996), 237-302.
\end{thebibliography}
\end{document}